\documentclass[]{article}

\usepackage{lipsum}
\usepackage[a4paper,margin=3.5cm]{geometry}
\usepackage{amsfonts}
\usepackage{graphicx}
\usepackage{epstopdf}
\usepackage{amsmath}          
\usepackage{xfrac}            
\usepackage{bm}               
\usepackage{booktabs}
\usepackage{color}
\usepackage{xcolor}
\usepackage{caption}
\usepackage{subcaption}
\usepackage{mathtools}
\usepackage{algorithm2e}
\usepackage[
  style=numeric,
  natbib,
  backend=biber]{biblatex}

\bibliography{./bibliography}

\providecommand{\keywords}[1]
{
  \def \and {$\; \cdot \;$}
  \\
  \small
  \noindent
  \textbf{Keywords } {#1}
}

\newenvironment{acknowledgements}
{
 \noindent
 \textbf{Acknowledgments }
}
{
  
}

\begin{document}
\title{A restricted additive smoother for finite cell flow problems \thanks{Financial support was provided by the German Research Foundation ({\it Deutsche Forschungsgemeinschaft, DFG}) in the framework of the collaborative research center SFB 837 {\it Interaction Modeling in Mechanized Tunneling}.}
}

\author{M. Saberi \thanks{High Performance Computing, Ruhr University Bochum, Universit\"{a}tsstr. 150, 44801 Bochum, Germany} \and
  A. Vogel \footnotemark[2]}

\maketitle

\begin{abstract}
In this work, we propose an adaptive geometric multigrid method for the solution of large-scale finite cell flow problems. The finite cell method seeks to circumvent the need for a boundary-conforming mesh through the embedding of the physical domain in a regular background mesh. As a result of the intersection between the physical domain and the background computational mesh, the resultant systems of equations are typically numerically ill-conditioned, rendering the appropriate treatment of cutcells a crucial aspect of the solver. To this end, we propose a smoother operator with favorable parallel properties and discuss its memory footprint and parallelization aspects. We propose three cache policies that offer a balance between cached and on-the-fly computation and discuss the optimization opportunities offered by the smoother operator. It is shown that the smoother operator, on account of its additive nature, can be replicated in parallel exactly with little communication overhead, which offers a major advantage in parallel settings as the geometric multigrid solver is consequently independent of the number of processes. The convergence and scalability of the geometric multigrid method is studied using numerical examples. It is shown that the iteration count of the solver remains bounded independent of the problem size and depth of the grid hierarchy. The solver is shown to obtain excellent weak and strong scaling using numerical benchmarks with more than 665 million degrees of freedom. The presented geometric multigrid solver is, therefore, an attractive option for the solution of large-scale finite cell problems in massively parallel high-performance computing environments.
\keywords{Finite cell method \and Stokes \and Geometric multigrid \and unfitted finite element method \and restricted additive}
\end{abstract}
\section{Introduction}
\label{sec:introduction}
The numerical approximation of partial differential equations using classical techniques such as the finite element method~\parencite{ern_2004}, the finite difference method~\parencite{leveque_2007} and the finite volume method~\parencite{eymard_2000} involves the generation of a boundary-conforming mesh, which is, especially for complex domains, a laborious and cost-intensive task. A conglomerate of non-conforming methods under the loose brand of unfitted finite element methods have recently emerged in response to such challenges and include techniques such as the extended finite element method~\parencite{belytschko_2001}, the finite cell method~\parencite{parvizian_2007,duster_2008,schillinger_2015}, the cut finite element method~\parencite{burman_2015} and the immersed finite element method~\parencite{peskin_2002,mittal_2005}. In this work, we focus on the finite cell method, where the physical domain of interest is embedded in a regular background computational mesh and recovered using a sufficiently accurate integration technique. The finite cell method has been successfully used for a variety of applications, including large deformation structural analysis~\parencite{schillinger_2012}, elasto-plasticity~\parencite{abedian_2013_finite,abedian_2014,ranjbar_2014}, linear and nonlinear fluid flow~\parencite{saberi_2023_adaptive}, thermo-elasticity~\parencite{zander_2012}, shell analysis~\parencite{rank_2012}, topology optimization~\parencite{parvizian_2012_topology}, etc., see also~\parencite{schillinger_2015} for a review. While the generation of a boundary conforming mesh is effectively circumvented in the finite cell method, and while the regular nature of the computational domain offers a wide range of optimization opportunities, the resultant system of equations is typically severely ill-conditioned mainly on account of unfavorably cut cells, rendering the solution of large-scale problems specially daunting as a consequence. Therefore, the development of efficient iterative solvers for the solution of the resultant system is a relevant challenge to the usage of such methods.

Geometric multigrid solvers are shown to be among the most efficient iterative methods for classical finite element methods~\parencite{hackbusch_1985} as well as the finite cell method~\parencite{deprenter_2019_multigrid,saberi_2020,saberi_2023_adaptive,jomo_2021}, capable of obtaining optimal convergence rates independent of the problem size. In this regard, the convergence of the geometric multigrid solver heavily depends on the effectiveness of the employed smoother operator, which is virtually always problem and discretization dependent. A tailored smoother operator for the treatment of cutcells was shown to be necessary for the convergence of geometric multigrid methods for the finite cell method, see, e.g.,~\parencite{deprenter_2019_multigrid,saberi_2020,saberi_2023_adaptive}.

In this work, we present an adaptive geometric multigrid method for the solution of the finite cell formulation of the Stokes problem. Furthermore, the saddle-point systems arising from the discretization of the Stokes equations appear also in a wide range of other domains such as power network analysis, finance and PDE optimization problems, making the presented geometric multigrid method for the solution of the Stokes equations relevant beyond flow applications. The Uzawa method~\parencite{elman_1994,bramble_1997_analysis}, multigrid methods~\parencite{verfurth_1984,hackbusch_1985,wittum_1989}, Krylov subspace solvers~\parencite{rusten_1992,silvester_1994,saad_2003} with a variety of preconditioning techniques including Schur complement methods~\parencite{patankar_1972,silvester_2001,kay_2002,elman_2006} and domain decomposition methods~\parencite{klawonn_1998_overlapping,klawonn_2000,pavarino_2000} are among the most popular solvers for the Stokes equations. A variety of smoothers have been proposed for the Stokes equations, including a smoother based on incomplete LU factorization~\parencite{wittum_1989}, the Braess-Sarazin smoother~\parencite{braess_1997} which is based on a global pressure correction scheme and the Vanka smoother~\parencite{vanka_1986} which is based on the solution of local saddle-point problems. The restricted additive Vanka (RAV) smoother was recently proposed in~\parencite{saberi_2022} and its parallel scalability for large-scale problems on distributed-memory computing clusters was studied in~\parencite{saberi_2023_efficient}. It was shown that the RAV smoother is competitive with the classical Vanka smoother in terms of convergence rate while being computationally less expensive; moreover, it can be represented exactly in parallel with minimal communication because of its additive nature. 

The aim of this work is to develop a scalable and efficient smoother for the finite cell formulation of the Stokes equations with focus on the solution of large-scale problems in high-performance computing environments. To this end, the main findings of~\parencite{saberi_2023_adaptive}, namely the appropriate treatment of cutcells in the finite cell method and the RAV smoother in~\parencite{saberi_2023_efficient} are used as inspiration for the development of the proposed smoother. The main contributions of this work can be summarized as follows:
\begin{itemize}
  \item A restricted additive smoother in the context of an adaptive geometric multigrid method for the finite cell formulation of the Stokes equations is proposed
  \item The computational and parallelization aspects of the smoother operator are discussed, and three cache policies are proposed. In particular, it is shown that the smoother operator can be replicated exactly in parallel with little communication overhead and offers memory optimization opportunities
  \item The parallel scalability and convergence of the geometric multigrid solver is studied using numerical benchmarks with up to more than 665 million degrees of freedom on distributed-memory machines
\end{itemize}
The remainder of this work is organized as follows. In Section~\ref{sec:fcm}, the finite cell formulation of the Stokes equations is presented. The geometric multigrid solver and the proposed smoother are discussed in detail in Section~\ref{sec:gmg}. The behavior of the geometric multigrid solver with the proposed smoother is studied in Section~\ref{sec:numerical_experiments} using a number of numerical experiments in terms of convergence, computational cost and weak and strong scalability in parallel on distributed-memory machines. Finally, conclusions are drawn in Section~\ref{sec:conclusions}.
\section{Finite cell formulation}
\label{sec:fcm}
In this section, we discuss the mixed finite cell formulation of the Stokes equations, whose strong form can be written as
\begin{equation}
\label{eq:stokes_strong}
\begin{split}
- \eta \nabla^{2}\bm{u} + \nabla p &= \bm{f} \quad \text{in} \; \Omega, \\
\nabla \cdot \bm{u} &= 0 \quad \text{in} \; \Omega,
\end{split}
\end{equation}
where $\Omega$ is the physical domain. $\bm{u}$ is the vector-valued fluid velocity and $p$ is the scalar fluid pressure. $\eta$ is the kinematic viscosity, and $\bm{f}$ is the body force function. Equation~\eqref{eq:stokes_strong} along with the boundary conditions
\begin{equation}
\label{eq:stokes_bc}
\begin{split}
\bm{u} &= \bm{w} \; \text{on} \; \Gamma_{D} \subset \partial \Omega, \\
\eta \frac{\partial \bm{u}}{\partial \bm{n}} - \bm{n} p &= \bm{h} \; \text{on} \; \Gamma_{N} := \partial \Omega \setminus \Gamma_{D}, \\
\end{split}
\end{equation}
forms the Stokes boundary-value problem, where $\partial \Omega$ is the boundary of the physical domain. $\Gamma_{D}$ and $\Gamma_{N}$ denote the Dirichlet and Neumann boundaries, respectively, such that $\partial \Omega = \Gamma_{D} \cup \Gamma_{N}$ and $\Gamma_{D} \cap \Gamma_{N} = \emptyset$. $\bm{w}$ and $\bm{h}$ are predefined functions. $\bm{n}$ is the outward normal vector to the boundary with unit length.

The boundary-conforming weak formulation of the Stokes equations can be obtained following the standard finite element procedure, see, e.g.,~\parencite{ern_2004} as:

Find $(\bm{u}, p) \in (\bm{V}_{w}, Q)$ such that
\begin{equation}
\begin{split}
\label{eq:stokes_weak_2}
(\eta \nabla \bm{v}, \nabla \bm{u})_{\Omega} - (\nabla \cdot \bm{v}, p)_{\Omega} - (q, \nabla \cdot \bm{u})_{\Omega} = \\
(\bm{v}, \bm{f})_{\Omega} + (\bm{v}, \bm{h})_{\Gamma_{N}}, \quad \forall (\bm{v}, q) \in (\bm{V}_{0}, Q), \\
\end{split}
\end{equation}
where $(\cdot, \cdot)_{\Omega}$ and $(\cdot, \cdot)_{\partial \Omega}$ denote the $L^{2}$ scalar product on $\Omega$ and $\partial \Omega$, respectively and
\begin{equation}
\begin{split}
\bm{V}_{w} &:= \{\bm{u} \in H^{1}(\Omega)^{d} \; | \; \bm{u} = \bm{w} \; \text{on} \; \Gamma_{D}\},\\
\bm{V}_{0} &:= \{\bm{v} \in H^{1}(\Omega)^{d} \; | \; \bm{v} = \bm{0} \; \text{on} \; \Gamma_{D}\},\\
Q &:= \{q \in L^{2}(\Omega) \},
\end{split}
\end{equation}
where $d$ is the spatial dimension of the domain.

We derive the weak finite cell formulation of the Stokes equations next, which involves on the one hand the extension of the physical domain to an embedding computational domain and the weak imposition of the essential boundary conditions on the other. The physical domain $\Omega$ is extended by a fictitious part to the embedding domain $\Omega_{e}$ in order to circumvent the need for a boundary-conforming mesh. $\Omega_{e}$ is typically chosen as regular such as to enable the straightforward generation and manipulation of the computational domain. We employ fully distributed space tree data structures, see, e.g.,~\parencite{burstedde_2011,bader_2012,saberi_2020} for the discretization of $\Omega_{e}$ that are capable of the efficient handling of regular spaces with low memory footprint and excellent load balancing. As a consequence of embedding the physical domain in $\Omega_{e}$, the boundary of the physical domain is not guaranteed to coincide with the computational domain, and essential boundary conditions must be imposed weakly in the finite cell formulation in contrast to the boundary-conforming case above, where the essential boundary conditions were included in the function spaces. Lagrange multipliers~\parencite{fernandez_2004,flemisch_2007,glowinski_2007,burman_2010}, the penalty method~\parencite{babuvska_1973,zhu_1998,burman_2010_ghost} and Nitsche's method~\parencite{nitsche_1971,hansbo_2002,dolbow_2009,embar_2010,burman_2012} are among the techniques used for the weak imposition of boundary conditions. We employ the symmetric Nitsche's method in this work. The weak finite cell formulation of the Stokes equations can be written as, see also~\parencite{saberi_2023_adaptive}:

Find $(\bm{u}, p) \in (\bm{V_{e}}, Q_{e})$ such that
\begin{equation}
\begin{split}
\label{eq:stokes_weak_fcm}
(\eta \nabla \bm{v}, \alpha \nabla \bm{u})_{\Omega_{e}} - (\bm{v}, \bm{n} \cdot \eta \nabla \bm{u})_{\Gamma_{D}} - (\bm{n} \cdot \eta \nabla \bm{v}, \bm{u})_{\Gamma_{D}} + (\bm{v}, \lambda \bm{u})_{\Gamma_{D}} \\
- (\nabla \cdot \bm{v}, \alpha p)_{\Omega_{e}} + (\bm{v}, \bm{n} p)_{\Gamma_{D}} + (\bm{n} q, \bm{u})_{\Gamma_{D}} - (q, \alpha \nabla \cdot \bm{u})_{\Omega_{e}} = \\
(\bm{v}, \alpha \bm{f})_{\Omega_{e}} + (\bm{v}, \bm{h})_{\Gamma_{N}} - (\bm{n} \cdot \eta \nabla \bm{v}, \bm{w})_{\Gamma_{D}} + (\bm{v}, \lambda \bm{w})_{\Gamma_{D}} + (\bm{n} q, \bm{w})_{\Gamma_{D}}, \\
\forall (\bm{v}, q) \in (\bm{V}_{e}, Q_{e}), \\
\end{split}
\end{equation}
where $\Omega_{e}$ is the embedding domain, $(\bm{V}_{e}, Q_{e})$ are given by
\begin{equation}
\begin{split}
\label{eq:space_unfitted_fem}
\bm{V}_{e} &:= \{\bm{u} \in H^{1}(\Omega_{e})^{d} \},\\
Q_{e} &:= \{q \in L^{2}(\Omega_{e}) \},
\end{split}
\end{equation}
and $\alpha$ is a scalar penalization parameter defined as
\begin{equation}
\begin{cases}
\alpha = 1 \quad \text{in} \; \Omega, \\
\alpha = 0 \quad \text{in} \; \Omega_{e} \setminus \Omega.
\end{cases}
\end{equation}
The penalization parameter $\alpha$ is typically chosen as a small value ($\alpha \ll 1$) instead in the fictitious part $\Omega_{e} \setminus \Omega$ in order to alleviate the severe numerical ill-conditioning of the resultant system of equations. $\lambda$ is the scalar stabilization parameter in Nitsche's method.

Introducing the discrete spaces $(\bm{V}_{e,h}, Q_{e,h}) \subset (\bm{V}_{e}, Q_{e})$, the bilinear and linear forms can be written as
\begin{equation}
\begin{split}
\label{eq:stokes_disc_weak_fcm_1}
a(\bm{v}_{h}, \bm{u}_{h}) + b(\bm{v}_{h}, p_{h}) &= f(\bm{v}_{h}) \\
b(q_{h}, \bm{u}_{h}) + c(q_{h}, p_{h}) &= g(q_{h}),\\
\end{split}
\end{equation}
where $a$, $b$, $f$ and $g$ are defined according to the weak form in Equation~\eqref{eq:stokes_weak_fcm} as
\begin{equation}
\begin{alignedat}{1}
\label{eq:stokes_disc_weak_fcm_2}
a(\bm{v}_{h}, \bm{u}_{h}) &:= (\eta \nabla \bm{v}_{h}, \alpha \nabla \bm{u}_{h})_{\Omega_{e, h}} - (\bm{v}_{h}, \bm{n} \cdot \eta \nabla \bm{u}_{h})_{\Gamma_{D, h}} - (\bm{n} \cdot \eta \nabla \bm{v}_{h}, \bm{u}_{h})_{\Gamma_{D, h}} \\
&+ (\bm{v}_{h}, \lambda \bm{u}_{h})_{\Gamma_{D, h}} \\
b(\bm{v}_{h}, p_{h}) &:= - (\nabla \cdot \bm{v}_{h}, \alpha p_{h})_{\Omega_{e, h}} + (\bm{v}_{h}, \bm{n} p_{h})_{\Gamma_{D, h}} \\
f(\bm{v}_{h}) &:= (\bm{v}_{h}, \alpha \bm{f})_{\Omega_{e, h}} + (\bm{v}_{h}, \bm{h})_{\Gamma_{N, h}} - (\bm{n} \cdot \eta \nabla \bm{v}_{h}, \bm{w})_{\Gamma_{D, h}} + (\bm{v}_{h}, \lambda \bm{w})_{\Gamma_{D, h}} \\
g(q_{h}) &:= (\bm{n} q_{h}, \bm{w})_{\Gamma_{D}}, \\
\end{alignedat}
\end{equation}
where $\Omega_{e, h}$ and $\Gamma_{\cdot, h}$ are appropriate discretizations of $\Omega_{e}$ and $\partial \Gamma$, respectively. An approximation of the embedding domain $\Omega_{e}$ is defined by $\Omega_{e, h}$ such that $\overline{\Omega}_{e,h} := \cup_{i=1}^{n_{K}} K_{i}$, where $\mathcal{T}_{e, h} := \{K_i\}_{i=1}^{n_K}$ with $\mathring{K}_i \cap \mathring{K}_j = \emptyset$ for $i \neq j$ form a tessellation of $\Omega_{e}$ into a set of $n_K$ compact, connected, Lipschitz sets $K_{i}$ with non-empty interior.

In this work, we focus on the $Q_{1}-Q_{1}$ discretization of the weak formulation. The equal-order $Q_{p}-Q_{p}$ elements do not satisfy the inf-sup condition, and an appropriate stabilization of the discrete form is, therefore, necessary. We note that equal order elements, especially of lower orders, nevertheless, remain attractive choices on account of their lower computational cost and ease of implementation compared to most of their conforming counterparts. The bilinear form $c$ in Equation~\eqref{eq:stokes_disc_weak_fcm_1} is then the stabilization term based on a local $L^{2}$ projection of the pressure~\parencite{dohrmann_2004,bochev_2006,li_2008}, given by
\begin{equation}
\label{eq:stabilization_l2}
c(q_{h}, p_{h}) := - \frac{1}{\eta}(q_{h} - \Pi_{0} q_{h}, p_{h} - \Pi_{0} p_{h})_{\Omega_{e, h}},
\end{equation}
where $\Pi_{0}$ is a local $L^{2}$-projection operator~\parencite{bochev_2006}, where the projection of a given function $q$ must satisfy
\begin{equation}
\int_{\Omega_{K}} (\Pi_{0}q - q) d \bm{x} = 0 \; \forall K \in \mathcal{T}_{e, h},
\end{equation}
where $\Omega_{K}$ is the domain associated with $K$.
\section{Geometric multigrid}
\label{sec:gmg}
In this section, we develop a monolithic adaptive geometric multigrid solver for the mixed finite cell formulation of Stokes equations on distributed-memory machines, where we briefly introduce the multigrid framework in the first part and focus on the proposed smoother operator in the second. The spatial discretization of the computational domain, as described in Section~\ref{sec:fcm}, is handled by space tree data structures, where adaptive mesh refinement towards the boundary of the physical domain and potentially other regions of interest is employed. Given a coarse cell, mesh refinement then corresponds to the division of the cell into its $2^{d}$ children, where $d$ is the spatial dimension. In addition, we impose a 2:1 balance such that neighboring cells are at most one level of refinement apart. The discrete computational domain may, therefore, include hanging nodes, which are handled as constraints and removed from the global system of equations. The discrete weak form in Equation~\eqref{eq:stokes_disc_weak_fcm_1} results in a linear system of the form
\begin{equation}
\label{eq:linear_system}
\mathbf{L} \mathbf{x} = \mathbf{b},
\end{equation}
where $\mathbf{L} := \begin{psmallmatrix}\mathbf{A} & \mathbf{B} \\ \mathbf{B}^{T} & \mathbf{C} \end{psmallmatrix}$ with $\mathbf{A} \in \mathbb{R}^{n_{u} \times n_{u}}$ and $\mathbf{B} \in \mathbb{R}^{n_{u} \times n_{p}}$ and $\mathbf{b} := \begin{psmallmatrix} \mathbf{f} \\ \mathbf{g} \end{psmallmatrix}$ with $\mathbf{f} \in \mathbb{R}^{n_{u}}$ and $\mathbf{g} \in \mathbb{R}^{n_{p}}$ are defined according to the bilinear and linear forms of the discrete weak problem in Equation~\eqref{eq:stokes_disc_weak_fcm_2}, respectively. The matrix $\mathbf{C} \in \mathbb{R}^{n_{p} \times n_{p}}$ is defined according to the stabilization term $c$ in Equation~\eqref{eq:stabilization_l2}. $\mathbf{x} := \begin{psmallmatrix} \mathbf{u} \\ \mathbf{p} \end{psmallmatrix}$ is the solution vector, where $\mathbf{u} \in \mathbb{R}^{n_{u}}$ and $\mathbf{p} \in \mathbb{R}^{n_{p}}$ are the coefficients of expansion of the velocity and pressure functions, respectively. Geometric multigrid methods seek the solution to the system of equations in Equation~\eqref{eq:linear_system} through a hierarchy of coarser discretizations, $\Omega_{e, h}^{l}, l = 1, \ldots, L$, where $\Omega_{e, h}^{L}$ and $\Omega_{e, h}^{1}$ are the finest and coarsest grids, respectively. We construct $\Omega_{e,h}^{l}, l = 1, \ldots, L$ top down, where given the grid $\Omega_{e,h}^{l}$, the immediate coarse grid $\Omega_{e,h}^{l-1}$ is obtained according to Algorithm~\ref{al:grid_hierarchy}. This process is repeated $L$ times, resulting in a nested grid hierarchy.

The multigrid cycle consists in the application of the smoother operator, transfer of the residual and correction vectors through the grid hierarchy and the solution of the coarsest problem. The grids $\Omega_{e,h}^{l}, l = 1, \ldots, L$ are fully distributed such that a given cell is uniquely owned by a single process. Therefore, it is desirable to keep the number of local cells per process roughly equal on each grid in order to maintain load balancing. As a consequence of adaptive mesh refinement, child cells and their corresponding parent are not guaranteed to remain on the same process after the application of load balancing. Therefore, data transfer between grids is carried out in two steps as follows. The restriction of a vector $\mathbf{v}^{l}$ on the fine grid $\Omega_{e,h}^{l}$ to $\mathbf{v}^{l - 1}$ on the coarse grid $\Omega_{e,h}^{l - 1}$, is achieved through the restriction operator $\mathcal{R}_{l}^{l - 1}: \Omega_{e,h}^{l} \to \Omega_{e,h}^{l - 1}$ which can be written as
\begin{equation}
\label{eq:restriction_op}
\mathcal{R}_{l}^{l - 1} := \mathcal{T}^{l - 1} \tilde{\mathcal{R}}_{l}^{l - 1}.
\end{equation}
$\mathcal{R}_{l}^{l - 1}$ is split into the intermediate restriction operator $\tilde{\mathcal{R}}_{l}^{l - 1}: \Omega_{e,h}^{l} \to \tilde{\Omega}_{e,h}^{l - 1}$ and the transfer operator $\mathcal{T}^{l - 1}: \tilde{\Omega}_{e,h}^{l - 1} \to \Omega_{e,h}^{l - 1}$, where $\tilde{\Omega}_{e,h}^{l - 1}$ is an intermediate process-local coarse grid, constructed from the fine grid $\Omega_{e,h}^{l}$ such as to guarantee that the parent of child nodes are on the same process as $\Omega_{e,h}^{l}$. Therefore, $\tilde{\mathcal{R}}_{l}^{l - 1}$ is a process-local operator that does not require any communication. The operator $\mathcal{T}^{l - 1}$ transfers $\tilde{\mathbf{v}}^{l - 1}$ to $\mathbf{v}^{l - 1}$ between $\tilde{\Omega}_{e,h}^{l - 1}$ and the distributed coarse grid $\Omega_{e,h}^{l - 1}$ and involves data communication. Prolongation of a vector $\mathbf{v}^{l - 1}$ on the coarse grid $\Omega_{e,h}^{l - 1}$ to $\mathbf{v}^{l}$ on the fine grid $\Omega_{e,h}^{l}$, expressed as $\mathcal{P}_{l - 1}^{l} := {\mathcal{R}_{l}^{l - 1}}^{T}$, similarly takes place in two steps, where $\mathbf{v}^{l - 1}$ is first transferred to an intermediate grid and subsequently to the fine grid.

The solution of the coarsest problem is the remaining step in the multigrid cycle, for which a direct solver is used in this work. We note that direct solvers typically suffer from high computational complexity on the one hand and suboptimal concurrency opportunities on the other; therefore, the efficient employment of direct solvers within the multigrid cycle in parallel presupposes that, through a sufficiently deep grid hierarchy, the coarsest problem is sufficiently small.

\begin{algorithm}[tb]
\SetKwInOut{Input}{input}\SetKwInOut{Output}{output}
\Input{$\Omega_{e,h}^{l}$}
\Output{$\Omega_{e,h}^{l-1}$}
\BlankLine
$\Omega_{e,h}^{l-1} \leftarrow \Omega_{e,h}^{l}$ \;
$r_{max} \leftarrow \text{max}(r_{K} \forall K \in \Omega_{e,h}^{l-1})$ \;
\BlankLine
\For{$K \in \Omega_{e,h}^{l-1}$}
{
    \If{$r_{K} == r_{max}$}
    {
        coarsen $K$ \;
    }
}
\BlankLine
Balance $\Omega_{e,h}^{l-1}$ \;
\BlankLine
\caption{Generation of the immediate coarse grid $\Omega_{e,h}^{l-1}$ from the fine grid $\Omega_{e,h}^{l}$. $r_{K}$ is the refinement level of $K$. Given a a cell $K$, the \texttt{coarsen} operator replaces $K$ and all of its siblings with their parent. The balance operator imposes the 2:1 balancing as explained in Section~\ref{sec:gmg}}
\label{al:grid_hierarchy}
\end{algorithm}
\subsection{The restricted additive smoother}
The convergence of multigrid methods heavily depends on the assumption that highly oscillatory frequencies of the error can be effectively eliminated on finer grids on the one hand and that slowly oscillatory frequencies can be adequately represented on coarser grids on the other. The smoother operator plays a major role in this regard, whose appropriate choice often determines the performance of the multigrid method.~\parencite{hackbusch_1985}

We propose a smoother for the solution of the finite cell formulation of Stokes equations in this section. We define the smoother operator first and discuss its parallelization aspects and memory footprint next, where three cache policies for the smoother operator are proposed.

A series of iterative corrections for the solution vector in Equation~\eqref{eq:linear_system} can be defined as
\begin{equation}
\label{eq:iteration}
\mathbf{x}^{k+1} = \mathbf{x}^{k} + \mathbf{S}(\mathbf{b} - \mathbf{L} \mathbf{x}^{k}),
\end{equation}
where $k$ denotes the iteration step and $\mathbf{S}$ denotes the smoother operator. The smoother operator can be written as
\begin{equation}
\label{eq:smoother}
\mathbf{S} = \sum_{i = 1}^{n_{p}} (\mathbf{\tilde R}_{i}^{T} \boldsymbol{\omega}_{i} \mathbf{L}_{i}^{-1}\mathbf{R}_{i}),
\end{equation}
and, in the context of Schwarz domain decomposition methods, see, e.g.,~\parencite{gander_2008}, consists in the application of corrections from a set of subdomains $\mathcal{S}_{i}, i = 1, \ldots, n_{p}$, where $n_{p}$ is the number of pressure nodes in the discrete domain. $\mathcal{S}_{i}$ is composed of the pressure degree of freedom $p_{i}$, all the pressure degrees of freedom that are connected to it, and all the velocity degrees of freedom that are either connected to $p_{i}$ or to the pressure degrees of freedom connected to $p_{i}$. The pressure DoF $p_{i}$ is said to be connected to a velocity DoF if the corresponding entry in the $i$-th row of $\mathbf{B}^{T}$ is non-zero. The subdomain restriction operator $\mathbf{R}_{i}$ effectively extracts the set of DoFs in $\mathcal{S}_{i}$, and is defined as those rows of the identity matrix $\mathbf{I} \in \mathbb{R}^{n \times n}$ that correspond to the DoFs in $\mathcal{S}_{i}$, where $n$ denotes the total number of degrees of freedom. $\mathbf{L}_{i} := \mathbf{R}_{i}\mathbf{L}\mathbf{R}_{i}^{T}$ is the subdomain block of $\mathbf{L}$ and $\boldsymbol{\omega}_{i}$ is a diagonal damping matrix. The prolongation operator $\mathbf{\tilde{R}}_{i}^{T}$ is an injection operator that extends a vector by padding it with zeros and only consists of the pressure DoF $p_{i}$ and the velocity DoFs on the same pressure node. The choice of the restriction and prolongation operators are inspired on the one hand by the observation that the treatment of cutcells in the finite cell formulation of the Stokes equations is essential for the convergence of multigrid methods~\parencite{saberi_2023_adaptive} and the favorable properties of the RAV smoother in~\parencite{saberi_2022} on the other.

The application of the correction from subdomain $\mathcal{S}_{i}$ consists in the solution of a local saddle-point problem, defined by the subdomain degrees of freedom, and corresponds to the local subdomain block $\mathbf{L}_{i}$. We note that in the general case, the local problem is dependent on material parameters on the one hand and the configuration of the elements in the local neighborhood on the other; therefore, we assume the subdomains to be spatially dependent.
\subsubsection{Parallelization aspects}
\label{sec:parallelization}
The proposed smoother is an additive operator, i.e., the subdomain corrections in Equation~\eqref{eq:smoother} are applied additively. In a parallel computing setting, the additive nature of the smoother operator is a major advantage as the subdomain corrections can effectively be applied simultaneously. The application of the subdomain corrections in parallel follows the ownership of the pressure degrees of freedom in the distributed computational domain, where $\mathcal{S}_{i}$ is owned by the process that owns $p_{i}$. Since pressure DoFs are uniquely owned, each subdomain is owned by exactly one process. Consequently, a given process is responsible for $n_{p}^{\texttt{proc}}$ subdomains, where $n_{p}^{\texttt{proc}}$ denotes the number of local pressure DoFs owned by the process.

In terms of data communication, we denote as process local a given subdomain if all of its degrees of freedom are owned by the same process and as off process otherwise. Process-local subdomains do not require any communication as all the necessary data for their application is available locally on the corresponding process. On the other hand, the degrees of freedom of off-process subdomains, which occur near process interfaces, are only partially owned by the corresponding process; therefore, off-process subdomains require the communication of the non-local degrees of freedom. We note that the pressure degrees of freedom in a given subdomain are either local or available through the ghost layer, and only velocity degrees of freedom belonging to subdomains on process interfaces may require data communication beyond the ghost layer.

Given the independence between subdomains, the smoother operator in Equation~\eqref{eq:smoother} can be replicated exactly in parallel, where the transfer of the non-local DoFs of off-process subdomains is the only necessary communication. Therefore, the smoother operator is independent of the number of processes in parallel. We note that $n_{\mathcal{S}}^{\texttt{offp}} \ll n_{p}^{\texttt{proc}}$, where $n_{\mathcal{S}}^{\texttt{offp}}$ is the number of off-process subdomains, and the associated communication overhead is, therefore, small.
\subsubsection{Cache policies}
\label{sec:cache}
The application of the subdomain corrections can be broken down into three steps: the retrieving of the local subdomain problem, the solution of the local problem and the application of the local correction. Given that the smoother operator in Equation~\eqref{eq:smoother} is typically applied several times, the optimization of these steps is clearly motivated from a computational perspective, leading naturally to three cache policies discussed below, see also,~\parencite{saberi_2023_efficient}. Given that the local subdomain problems are relatively small, we solve the local problems down to machine accuracy by computing the LU factorization of the local block $\mathbf{L}_{i}$ using a direct solver.

The first of the three cache policies, denoted as \texttt{cache\_matrix}, concerns itself with the elimination of the first step above, where the local subdomain problem is retrieved and stored once during the initialization of the smoother operator. Assuming that the global system of equations is stored in a sparse matrix format, the retrieving of the subdomain degrees of freedom, which are not, in general, stored in adjacent rows, is often costly since sparse matrix formats are typically not optimized for such operations. We further note that the local subdomain problem is itself sparse and can, therefore, be stored using a sparse format in order to minimize the memory footprint of the \texttt{cache\_matrix} policy. The application of a given subdomain correction using the \texttt{cache\_matrix} policy then consists in the inversion of the local block and updating the local residual vector with the local subdomain correction.

A more aggressive approach, denoted as \texttt{cache\_inverse}, is to compute and store the inverse of the local subdomain problems during the initialization of the smoother operator. The inversion of the local blocks $\mathbf{L}_{i}$ is by far the most computationally intensive step in the application of the subdomain corrections; therefore, the \texttt{cache\_inverse} policy is expected to result in a significant computational gain. On the other hand, it should be noted that the inverse of the local subdomain problem is in general a dense matrix, and hence at first glance, the \texttt{cache\_inverse} policy may seem to entail a heavy penalty in terms of memory consumption. However, the prolongation operator $\tilde{\mathbf{R}}^{T}_{i}$ of the proposed smoother operator offers an optimization opportunity upon closer inspection, namely as the prolongation operator of subdomain $\mathcal{S}_{i}$ contains only the pressure DoF $p_{i}$ and the velocity DoFs on the same node, it is only necessary to store the corresponding rows from the local inverse. In other words, the \texttt{cache\_inverse} policy requires the storage of $1 + d$ rows of the inverse of the local subdomain problem, where $d$ is the spatial dimension.

The cache policies above favor computational efficiency at the cost of increasing the memory footprint of the smoother operator. However, it may, in particular cases, be desirable to reduce memory consumption as far as possible, e.g., when the available main memory is extremely constrained. The corresponding cache policy is denoted as \texttt{cache\_none}, where all the steps associated with the application of the smoother operator are performed on the fly without the need to cache any extra matrices. Given the computational cost associated with the retrieving and inversion of the local subdomain problem, the presented cache policies are expected to offer a balance between computational efficiency and memory footprint. The cache policies are further investigated in Section~\ref{sec:numerical_experiments}.
\section{Numerical experiments}
\label{sec:numerical_experiments}
\begin{figure}[t]
\centering
\includegraphics{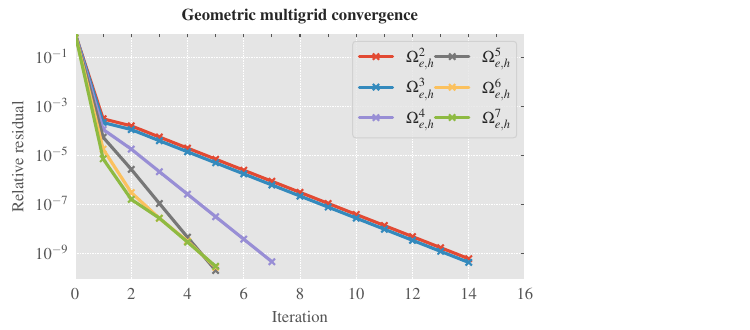}
\caption{The convergence of the geometric multigrid solver in the mesh study, where the grid hierarchy in Table~\ref{tab:grid} is solved. $\Omega_{e,h}^{l}, l = 2 \ldots 7$ denotes the fine problem. All problems employ $\Omega_{e,h}^{1}$ as the coarse grid; therefore, finer problems use a deeper grid hierarchy. The reduction of the relative residual by $10^{9}$ is used as the convergence criterion}
\label{fig:gmg_convergence}
\end{figure}
\begin{figure}[t]
\centering
\includegraphics{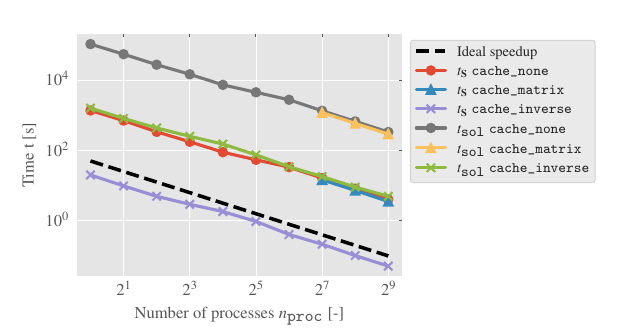}
\caption{The strong scaling of the geometric multigrid solver using different cache policies in the channel flow benchmark with $\Omega_{e,h}^{5}$ as the fine grid, see Table~\ref{tab:grid} with $n_{\texttt{proc}} = 1, \ldots, 512$. $t_{\mathbf{S}}$ denotes the runtime of the smoother per iteration on the fine grid, and $t_{\text{sol}}$ denotes the total solver runtime including the setup time. We note that the required memory by the \texttt{cache\_matrix} policy exceeds the available main memory of up to two compute nodes}
\label{fig:strong_scaling}
\end{figure}
\begin{figure}[t]
\centering
\includegraphics{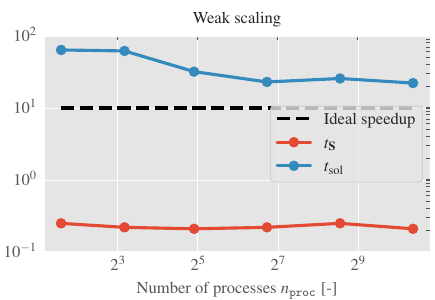}
\caption{The weak scaling of the geometric multigrid solver in the channel flow benchmark using the grid hierarchy in Table~\ref{tab:grid}, where the number of degrees of freedom per process is kept roughly constant. $t_{\mathbf{S}}$ denotes the runtime of the smoother per iteration on the fine grid, and $t_{\text{sol}}$ denotes the total solver runtime}
\label{fig:weak_scaling}
\end{figure}
\begin{table}[t]
\caption{The grid hierarchy of the channel flow benchmark, used for the grid study as well as the strong and weak scaling studies. $n_{c}$ and $n_{\text{DoF}}$ denote the number of cells and degrees of freedom, respectively. $n_{\text{proc}}$ is the number of processes used in the weak scaling study and $n_{\text{it}}$ is the iteration count of the solver}
\label{tab:grid}
    \centering
    \begin{tabular}[b]{|l|r|r|r|r|r|}
    \toprule
    Grid & $n_{c}$ & $n_{\text{DoF}}$ & $\sfrac{n_{\text{DoF}}^{l}}{n_{\text{DoF}}^{l - 1}}$ & $n_{\text{proc}}$ & $n_{\text{it}}$ \\
    \midrule
    $\Omega_{e,h}^{1}$ & 262,144 & 790,275 & - & - & -\\
    $\Omega_{e,h}^{2}$ & 526,516 & 1,583,127 & 2.00 & 3 & 14\\
    $\Omega_{e,h}^{3}$ & 1,491,076 & 4,477,533 & 2.83 & 9 & 14\\
    $\Omega_{e,h}^{4}$ & 4,999,612 & 15,004,455 & 3.35 & 30 & 7\\
    $\Omega_{e,h}^{5}$ & 17,669,596 & 53,016,765 & 3.53 & 106 & 5\\
    $\Omega_{e,h}^{6}$ & 62,884,444 & 188,665,497 & 3.56 & 376 & 5\\
    $\Omega_{e,h}^{7}$ & 221,782,588 & 665,367,261 & 3.53 & 1330 & 5\\
    \bottomrule
    \end{tabular}
\end{table}
In this section, we study the performance and scalability of the presented geometric multigrid solver using a number of numerical benchmarks. More specifically, the convergence of the geometric multigrid method is examined using a mesh study in the first part. The proposed cache policies in Section~\ref{sec:cache} are studied next in terms of computational efficiency and memory footprint, and the weak and strong scaling of the geometric multigrid solver for large-scale problems is presented on a distributed-memory cluster.

We employ as benchmark the well-known channel flow problem, see~\parencite{schafer_1996}, where the inflow with a prescribed velocity profile enters a channel with dimensions $2.2 m \times 0.41 m$ with a cylindrical obstacle with a diameter of $0.1 m$ near the inflow region. A no-slip wall condition is imposed on the bottom and top sides of the channel. The kinematic viscosity of the fluid is $\eta = \frac{1}{1000}\;\sfrac{m^{2}}{s}$. The inflow profile is described as
\begin{equation}
\label{eq:cf_inflow}
\begin{alignedat}{3}
u_{x}(0, y) &= \frac{4 \bar{u} y (H - y)}{H^{2}} \quad &&y \in [0, H], \\
u_{y}(0, y) &= 0 \quad &&y \in [0, H],
\end{alignedat}
\end{equation}
where $\bar{u} = 0.3 \; \sfrac{m}{s}$ and $H$ is the width of the channel. The penalization parameter $\alpha$ in the finite cell method is chosen as $10^{-10}$.

The computational grids, as shown in Table~\ref{tab:grid}, are produced using adaptive mesh refinement towards the boundary layers, i.e., the inflow and outflow regions, the top and bottom walls and the surface of the cylinder. Geometric multigrid is used as a standalone solver in order to better isolate its performance. The reduction of the relative residual by $10^{9}$ is used as the convergence criterion. We note that multigrid methods can be used as preconditioners in Krylov subspace methods, which often leads to improved convergence. A multigrid V cycle with six steps of pre- and post-smoothing is employed. A damping factor of 0.8, which was the minimum amount observed to be necessary, is applied to the smoother. We note that fewer smoothing steps were sufficient for convergence; nevertheless, more smoothing steps led to more rapid convergence and, given the relatively low computational cost of the smoother operator, moderately more smoothing steps were observed to provide the solver with a cost effective improvement. The reported runtime per smoother iteration consists of the application of the smoother operator as well as the synchronization of the residual vector.

The numerical experiments are carried out on a distributed-memory CPU cluster, where each node is equipped with double-socket Intel Xeon Skylake Gold 6148 CPUs with 20 cores running at 2.4 GHz and 180 GB of DDR4 main memory. The computing nodes are connected via a 100 GBit/s Intel Omni-Path Interconnect. The numerical examples are run in pure MPI mode, where each node is filled with up to 32 processes. An in-house C++ implementation is employed for the numerical experiments in the following. p4est~\parencite{burstedde_2011} and PETSc~\parencite{petsc_1997} are used for grid manipulation and some linear algebra functionalities, respectively.

We focus on the convergence of the geometric multigrid solver first, where a series of progressively finer problems with deeper grid hierarchies are solved. The grid hierarchy in Table~\ref{tab:grid} is employed, where the fine grids $\Omega_{e,h}^{l}, l = 2 \ldots 7$ are solved using $\Omega_{e,h}^{1}$ as the coarse grid; therefore, larger problems employ a deeper grid hierarchy. The convergence of the GMG solver is shown in Figure~\ref{fig:gmg_convergence}. It can be seen that the iteration count of the solver remains bounded independent of both the problem size and the depth of the grid hierarchy. Nevertheless, the solver tends to favor deeper grid hierarchies, and the iteration count is observed to improve for finer problems, which employ a deeper hierarchy, until it is stabilized. We note that aside from the depth of the grid hierarchy, the observed behavior is associated with the improved resolution of the boundary layers on finer grids.

We focus on the scalability of the geometric multigrid solver and the performance of the cache policies in Section~\ref{sec:cache} next. The strong scaling of the solver with $\Omega_{e,h}^{5}$ as the fine grid is shown in Figure~\ref{fig:strong_scaling}, where the runtime per smoother iteration as well as the total solver runtime using an increasing number of processes is reported. The performance of the solver with the different cache policies is shown in the same figure, where it can be observed that \texttt{cache\_inverse} is by far the most efficient policy in terms of total computational runtime. The \texttt{cache\_matrix} policy is slightly more efficient than the \texttt{cache\_none} policy; however, at a significantly higher cost in terms of memory footprint, e.g., the memory requirement of \texttt{cache\_matrix} exceeds the available main memory of up to two compute nodes. Finally, \texttt{cache\_none} does not require the storage of any extra matrices, and therefore, offers an alternative with the lowest memory footprint; however, as the retrieving and solution of the subdomain problems are performed on the fly, its computational cost is significantly higher in comparison. We note that in the absence of the memory optimization for the \texttt{cache\_inverse} policy, explained in Section~\ref{sec:cache}, the memory footprint of \texttt{cache\_inverse} would be even larger than \texttt{cache\_matrix}. However, given this memory optimization opportunity and given the massive computational gain compared to \texttt{cache\_none} and \texttt{cache\_matrix} policies, the \texttt{cache\_inverse} should clearly be preferred for the presented smoother operator minus exceptional cases where the available main memory in extremely limited.

It can be observed in Figure~\ref{fig:strong_scaling} that both the runtime of the smoother operator and the total solver runtime closely follow the ideal speedup. We note that as the smoother operator is exactly replicated in parallel, see Section~\ref{sec:parallelization}, the convergence and iteration count of the solver is independent of the number of processes.

The weak scaling of the solver with the \texttt{cache\_inverse} policy is presented is Figure~\ref{fig:weak_scaling}, where the grid hierarchy in Table~\ref{tab:grid} is solved using a progressively increasing number of processes such that the number of degrees of freedom per process is roughly constant. For a given number of processes, the smallest number of nodes is used, where each node is filled with up to 32 MPI processes. The runtime of the smoother operator is observed to remain virtually constant, and the total solver runtime remains roughly constant. The observed difference in the total solver runtime between the three coarsest problems and the finer problems can be imputed mainly to the difference in the iteration count of the solver, see Table~\ref{tab:grid}. In addition, we note that the large variation between the fine and coarse grids may lead to micro-parallelization on the coarser grids, where there is an imbalance between the number of operations and the number of utilized processes, and may, therefore, affect the scalability of multigrid methods. In such cases, an agglomeration technique, where fewer processes are assigned to coarser grids, is typically employed. Nevertheless, the geometric multigrid solver demonstrates excellent weak scaling for the fine problems in the presented benchmarks.
\section{Conclusions}
\label{sec:conclusions}
In this work, we proposed an adaptive geometric multigrid method for the solution of the finite cell formulation of the Stokes equations. The finite cell method, mainly on account of the existence of cutcells where the physical domain intersects the background computational mesh, leads to systems of equations that are numerically ill-conditioned; therefore, the appropriate treatment of cutcells is a crucial aspect of the multigrid method, see also~\parencite{saberi_2023_adaptive}. We proposed a smoother with favorable properties from a parallel computing perspective, and discussed its computational and parallelization aspects. The smoother operator, on account of its additive nature, can be replicated in parallel exactly with little communication overhead; therefore, the GMG solver with the proposed smoother is independent of the number of processes. Furthermore, three cache policies for the smoother operator were proposed, offering a balance between efficiency and on-the-fly computation. It is shown that the \texttt{cache\_inverse} policy, where the inverse of the local subdomain problems are computed and stored once during the initialization of the smoother operator, is by far the most computationally efficient policy in comparison. Moreover, the smoother operator is shown to offer an optimization opportunity in terms of memory footprint for the \texttt{cache\_inverse} policy; therefore, the \texttt{cache\_inverse} policy is clearly preferred except in rare cases where the available main memory is extremely constrained. The convergence and parallel scalability of the geometric multigrid solver with the presented smoother is studied using numerical benchmarks on a distributed-memory cluster. The convergence of the GMG solver is shown to remain bounded independent of both the problem size and the depth of the grid hierarchy. The strong and weak scaling of the solver is observed to closely follow the ideal speedup using numerical examples with more than 665 million degrees of freedom. Given the small communication overhead for the exact replication of the smoother operator in parallel and the excellent efficiency and scalability of the solver, the proposed geometric multigrid solver is an attractive option for the solution of large-scale finite cell problems on high-performance computing machines.

\begin{acknowledgements}
Financial support was provided by the German Research Foundation ({\it Deutsche Forschungsgemeinschaft, DFG}) in the framework of subproject C4 of the Collaborative Research Center SFB 837 {\it Interaction Modeling in Mechanized Tunneling}, grant number 77309832. This support is gratefully acknowledged. We also gratefully acknowledge the computing time on the computing cluster of the SFB 837 and the Department of Civil and Environmental Engineering at Ruhr University Bochum, which was used for the presented numerical studies.
\end{acknowledgements}

\printbibliography

@Article{eymard_2000,
  author    = {Eymard, Robert and Gallou{\"e}t, Thierry and Herbin, Rapha{\`e}le},
  title     = {Finite volume methods},
  pages     = {713--1018},
  volume    = {7},
  journal   = {Handbook of numerical analysis},
  publisher = {Elsevier},
  year      = {2000},
}

@Article{silvester_2001,
  author    = {Silvester, David and Elman, Howard and Kay, David and Wathen, Andrew},
  title     = {Efficient preconditioning of the linearized Navier--Stokes equations for incompressible flow},
  number    = {1-2},
  pages     = {261--279},
  volume    = {128},
  journal   = {Journal of Computational and Applied Mathematics},
  publisher = {Elsevier},
  year      = {2001},
}

@InCollection{petsc_1997,
  author    = {Satish Balay and William D. Gropp and Lois Curfman McInnes and Barry F. Smith},
  booktitle = {Modern Software Tools in Scientific Computing},
  title     = {Efficient Management of Parallelism in Object Oriented Numerical Software Libraries},
  editor    = {E. Arge and A. M. Bruaset and H. P. Langtangen},
  pages     = {163--202},
  publisher = {Birkh{\"{a}}user Press},
  year      = {1997},
}

@Article{abedian_2013_finite,
  author    = {Abedian, Alireza and Parvizian, Jamshid and D{\"u}ster, Alexander and Rank, Ernst},
  title     = {The finite cell method for the J2 flow theory of plasticity},
  pages     = {37--47},
  volume    = {69},
  journal   = {Finite Elements in Analysis and Design},
  publisher = {Elsevier},
  year      = {2013},
}

@Book{ern_2004,
  author    = {Ern, Alexandre and Guermond, Jean-Luc},
  title     = {Theory and practice of finite elements},
  publisher = {Springer},
  volume    = {159},
  year      = {2004},
}

@Article{babuvska_1973,
  author   = {Babu{\v{s}}ka, Ivo},
  title    = {The finite element method with penalty},
  number   = {122},
  pages    = {221--228},
  volume   = {27},
  journal  = {Mathematics of computation},
  year     = {1973},
}

@Article{silvester_1994,
  author    = {Silvester, David and Wathen, Andrew},
  title     = {Fast iterative solution of stabilised Stokes systems part II: using general block preconditioners},
  number    = {5},
  pages     = {1352--1367},
  volume    = {31},
  journal   = {SIAM Journal on Numerical Analysis},
  publisher = {SIAM},
  year      = {1994},
}

@Book{saad_2003,
  author    = {Saad, Yousef},
  title     = {Iterative methods for sparse linear systems},
  publisher = {siam},
  volume    = {82},
  year      = {2003},
}

@Article{klawonn_2000,
  author    = {Klawonn, Axel and Pavarino, Luca F},
  title     = {A comparison of overlapping Schwarz methods and block preconditioners for saddle point problems},
  number    = {1},
  pages     = {1--25},
  volume    = {7},
  journal   = {Numerical linear algebra with applications},
  publisher = {Wiley Online Library},
  year      = {2000},
}

@Article{schillinger_2012,
  author    = {Schillinger, Dominik and Ruess, Martin and Zander, Nils and Bazilevs, Yuri and D{\"u}ster, Alexander and Rank, Ernst},
  title     = {Small and large deformation analysis with the p-and B-spline versions of the Finite Cell Method},
  number    = {4},
  pages     = {445--478},
  volume    = {50},
  journal   = {Computational Mechanics},
  publisher = {Springer},
  year      = {2012},
}

@Article{schillinger_2015,
  author    = {Schillinger, Dominik and Ruess, Martin},
  title     = {The Finite Cell Method: A review in the context of higher-order structural analysis of CAD and image-based geometric models},
  number    = {3},
  pages     = {391--455},
  volume    = {22},
  journal   = {Archives of Computational Methods in Engineering},
  publisher = {Springer},
  year      = {2015},
}

@Article{burman_2012,
  author    = {Burman, Erik and Hansbo, Peter},
  title     = {Fictitious domain finite element methods using cut elements: II. A stabilized Nitsche method},
  number    = {4},
  pages     = {328--341},
  volume    = {62},
  journal   = {Applied Numerical Mathematics},
  publisher = {Elsevier},
  year      = {2012},
}

@Article{abedian_2014,
  author    = {Abedian, Alireza and Parvizian, Jamshid and D{\"u}ster, Alexander and Rank, Ernst},
  title     = {Finite cell method compared to h-version finite element method for elasto-plastic problems},
  pages     = {1239--1248},
  volume    = {35},
  journal   = {Applied Mathematics and Mechanics},
  publisher = {Springer},
  year      = {2014},
}

@Article{burman_2010,
  author    = {Burman, Erik and Hansbo, Peter},
  title     = {Fictitious domain finite element methods using cut elements: I. A stabilized Lagrange multiplier method},
  number    = {41-44},
  pages     = {2680--2686},
  volume    = {199},
  journal   = {Computer Methods in Applied Mechanics and Engineering},
  publisher = {Elsevier},
  year      = {2010},
}

@Article{burman_2015,
  author    = {Burman, Erik and Claus, Susanne and Hansbo, Peter and Larson, Mats G and Massing, Andr{\'e}},
  title     = {CutFEM: discretizing geometry and partial differential equations},
  number    = {7},
  pages     = {472--501},
  volume    = {104},
  journal   = {International Journal for Numerical Methods in Engineering},
  publisher = {Wiley Online Library},
  year      = {2015},
}

@Article{hansbo_2002,
  author    = {Hansbo, Anita and Hansbo, Peter},
  title     = {An unfitted finite element method, based on Nitsche’s method, for elliptic interface problems},
  number    = {47-48},
  pages     = {5537--5552},
  volume    = {191},
  comment   = {They only address interface problems. Dirichlet boundary conditions are not considered. The Nitsche's method is used.},
  journal   = {Computer methods in applied mechanics and engineering},
  publisher = {Elsevier},
  year      = {2002},
}

@Article{saberi_2023_adaptive,
  author    = {Saberi, S. and Meschke, G. and Vogel, A.},
  title     = {Adaptive geometric multigrid for the mixed finite cell formulation of Stokes and Navier–Stokes equations},
  doi       = {https://doi.org/10.1002/fld.5180},
  eprint    = {https://onlinelibrary.wiley.com/doi/pdf/10.1002/fld.5180},
  number    = {7},
  pages     = {1035--1053},
  url       = {https://onlinelibrary.wiley.com/doi/abs/10.1002/fld.5180},
  volume    = {95},
  journal   = {International Journal for Numerical Methods in Fluids},
  publisher = {Wiley Online Library},
  year      = {2023},
}

@Book{leveque_2007,
  author    = {LeVeque, Randall J},
  title     = {Finite difference methods for ordinary and partial differential equations: steady-state and time-dependent problems},
  publisher = {SIAM},
  year      = {2007},
}

@Article{dohrmann_2004,
  author    = {Dohrmann, Clark R and Bochev, Pavel B},
  title     = {A stabilized finite element method for the Stokes problem based on polynomial pressure projections},
  number    = {2},
  pages     = {183--201},
  volume    = {46},
  journal   = {International Journal for Numerical Methods in Fluids},
  publisher = {Wiley Online Library},
  year      = {2004},
}

@InCollection{patankar_1972,
  author    = {Patankar, Suhas V and Spalding, D Brian},
  booktitle = {Numerical prediction of flow, heat transfer, turbulence and combustion},
  title     = {A calculation procedure for heat, mass and momentum transfer in three-dimensional parabolic flows},
  pages     = {54--73},
  publisher = {Elsevier},
  year      = {1972},
}

@InProceedings{saberi_2023_efficient,
  author       = {Saberi, S and Meschke, G and Vogel, A},
  booktitle    = {European Conference on Parallel Processing},
  title        = {An Efficient Parallel Adaptive GMG Solver for Large-Scale Stokes Problems},
  organization = {Springer},
  pages        = {694--709},
  year         = {2023},
}

@Article{fernandez_2004,
  author    = {Fern{\'a}ndez-M{\'e}ndez, Sonia and Huerta, Antonio},
  title     = {Imposing essential boundary conditions in mesh-free methods},
  number    = {12-14},
  pages     = {1257--1275},
  volume    = {193},
  journal   = {Computer methods in applied mechanics and engineering},
  publisher = {Elsevier},
  year      = {2004},
}

@Article{zander_2012,
  author    = {Zander, Nils and Kollmannsberger, Stefan and Ruess, Martin and Yosibash, Zohar and Rank, Ernst},
  title     = {The finite cell method for linear thermoelasticity},
  number    = {11},
  pages     = {3527--3541},
  volume    = {64},
  journal   = {Computers \& Mathematics with Applications},
  publisher = {Elsevier},
  year      = {2012},
}

@Article{dolbow_2009,
  author    = {Dolbow, John and Harari, Isaac},
  title     = {An efficient finite element method for embedded interface problems},
  number    = {2},
  pages     = {229--252},
  volume    = {78},
  journal   = {International journal for numerical methods in engineering},
  publisher = {Wiley Online Library},
  year      = {2009},
}

@Article{jomo_2021,
  author    = {Jomo, John and Oztoprak, Oguz and de Prenter, Frits and Zander, Nils and Kollmannsberger, Stefan and Rank, Ernst},
  title     = {Hierarchical multigrid approaches for the finite cell method on uniform and multi-level hp-refined grids},
  pages     = {114075},
  volume    = {386},
  journal   = {Computer Methods in Applied Mechanics and Engineering},
  publisher = {Elsevier},
  year      = {2021},
}

@Article{zhu_1998,
  author    = {Zhu, T and Atluri, SN},
  title     = {A modified collocation method and a penalty formulation for enforcing the essential boundary conditions in the element free Galerkin method},
  number    = {3},
  pages     = {211--222},
  volume    = {21},
  journal   = {Computational Mechanics},
  publisher = {Springer},
  year      = {1998},
}

@Article{bramble_1997_analysis,
  author    = {Bramble, James H and Pasciak, Joseph E and Vassilev, Apostol T},
  title     = {Analysis of the inexact Uzawa algorithm for saddle point problems},
  number    = {3},
  pages     = {1072--1092},
  volume    = {34},
  journal   = {SIAM Journal on Numerical Analysis},
  publisher = {SIAM},
  year      = {1997},
}

@Article{verfurth_1984,
  author    = {Verf{\"u}rth, R{\"u}diger},
  title     = {A multilevel algorithm for mixed problems},
  number    = {2},
  pages     = {264--271},
  volume    = {21},
  journal   = {SIAM journal on numerical analysis},
  publisher = {SIAM},
  year      = {1984},
}

@Article{rank_2012,
  author    = {Rank, Ernst and Ruess, Martin and Kollmannsberger, Stefan and Schillinger, Dominik and D{\"u}ster, Alexander},
  title     = {Geometric modeling, isogeometric analysis and the finite cell method},
  pages     = {104--115},
  volume    = {249},
  journal   = {Computer Methods in Applied Mechanics and Engineering},
  publisher = {Elsevier},
  year      = {2012},
}

@Article{peskin_2002,
  author    = {Peskin, Charles S},
  title     = {The immersed boundary method},
  pages     = {479--517},
  volume    = {11},
  journal   = {Acta numerica},
  publisher = {Cambridge University Press},
  year      = {2002},
}

@Article{parvizian_2012_topology,
  author    = {Parvizian, J and D{\"u}ster, A and Rank, E},
  title     = {Topology optimization using the finite cell method},
  number    = {1},
  pages     = {57--78},
  volume    = {13},
  journal   = {Optimization and Engineering},
  publisher = {Springer},
  year      = {2012},
}

@Article{rusten_1992,
  author    = {Rusten, Torgeir and Winther, Ragnar},
  title     = {A preconditioned iterative method for saddlepoint problems},
  number    = {3},
  pages     = {887--904},
  volume    = {13},
  journal   = {SIAM Journal on Matrix Analysis and Applications},
  publisher = {SIAM},
  year      = {1992},
}

@Article{bochev_2006,
  author    = {Bochev, Pavel B and Dohrmann, Clark R and Gunzburger, Max D},
  title     = {Stabilization of low-order mixed finite elements for the Stokes equations},
  number    = {1},
  pages     = {82--101},
  volume    = {44},
  journal   = {SIAM Journal on Numerical Analysis},
  publisher = {SIAM},
  year      = {2006},
}

@Article{kay_2002,
  author    = {Kay, David and Loghin, Daniel and Wathen, Andrew},
  title     = {A preconditioner for the steady-state Navier--Stokes equations},
  number    = {1},
  pages     = {237--256},
  volume    = {24},
  journal   = {SIAM Journal on Scientific Computing},
  publisher = {SIAM},
  year      = {2002},
}

@Book{bader_2012,
  author    = {Bader, Michael},
  title     = {Space-filling curves: an introduction with applications in scientific computing},
  publisher = {Springer Science \& Business Media},
  volume    = {9},
  year      = {2012},
}

@Article{belytschko_2001,
  author    = {Belytschko, Ted and Mo{\"e}s, Nicolas and Usui, Shuji and Parimi, Chandu},
  title     = {Arbitrary discontinuities in finite elements},
  number    = {4},
  pages     = {993--1013},
  volume    = {50},
  journal   = {International Journal for Numerical Methods in Engineering},
  publisher = {Wiley Online Library},
  year      = {2001},
}

@Article{burstedde_2011,
  author    = {Burstedde, Carsten and Wilcox, Lucas C and Ghattas, Omar},
  title     = {p4est: Scalable algorithms for parallel adaptive mesh refinement on forests of octrees},
  number    = {3},
  pages     = {1103--1133},
  volume    = {33},
  journal   = {SIAM Journal on Scientific Computing},
  publisher = {SIAM},
  year      = {2011},
}

@Book{hackbusch_1985,
  author    = {Hackbusch, Wolfgang},
  title     = {Multi-grid methods and applications},
  edition   = {1},
  publisher = {Springer Berlin, Heidelberg},
  series    = {Springer Series in Computational Mathematics},
  year      = {1985},
}

@Article{ranjbar_2014,
  author    = {Ranjbar, M and Mashayekhi, M and Parvizian, J and D{\"u}ster, A and Rank, E},
  title     = {Using the finite cell method to predict crack initiation in ductile materials},
  pages     = {427--434},
  volume    = {82},
  journal   = {Computational Materials Science},
  publisher = {Elsevier},
  year      = {2014},
}

@Article{saberi_2022,
  author    = {Saberi, S. and Meschke, G. and Vogel, A.},
  title     = {A restricted additive Vanka smoother for geometric multigrid},
  pages     = {111123},
  volume    = {459},
  journal   = {Journal of Computational Physics},
  publisher = {Elsevier},
  year      = {2022},
}

@Article{elman_2006,
  author    = {Elman, Howard and Howle, Victoria E and Shadid, John and Shuttleworth, Robert and Tuminaro, Ray},
  title     = {Block preconditioners based on approximate commutators},
  number    = {5},
  pages     = {1651--1668},
  volume    = {27},
  journal   = {SIAM Journal on Scientific Computing},
  publisher = {SIAM},
  year      = {2006},
}

@InProceedings{saberi_2020,
  author       = {Saberi, S. and Vogel, A. and Meschke, G.},
  booktitle    = {Euro-Par 2020: Parallel Processing},
  title        = {Parallel Finite Cell Method with Adaptive Geometric Multigrid},
  editor       = {Malawski, Maciej and Rzadca, Krzysztof},
  isbn         = {978-3-030-57675-2},
  organization = {Springer},
  pages        = {578--593},
  publisher    = {Springer International Publishing},
  address      = {Cham},
  year         = {2020},
}

@Article{braess_1997,
  author    = {Braess, Dietrich and Sarazin, Regina},
  title     = {An efficient smoother for the Stokes problem},
  number    = {1},
  pages     = {3--19},
  volume    = {23},
  journal   = {Applied Numerical Mathematics},
  publisher = {Elsevier},
  year      = {1997},
}

@InCollection{schafer_1996,
  author    = {Sch{\"a}fer, Michael and Turek, Stefan and Durst, Franz and Krause, Egon and Rannacher, Rolf},
  booktitle = {Flow simulation with high-performance computers II},
  title     = {Benchmark computations of laminar flow around a cylinder},
  pages     = {547--566},
  publisher = {Springer},
  year      = {1996},
}

@Article{elman_1994,
  author    = {Elman, Howard C and Golub, Gene H},
  title     = {Inexact and preconditioned Uzawa algorithms for saddle point problems},
  number    = {6},
  pages     = {1645--1661},
  volume    = {31},
  journal   = {SIAM Journal on Numerical Analysis},
  publisher = {SIAM},
  year      = {1994},
}

@InProceedings{nitsche_1971,
  author       = {Nitsche, Joachim},
  booktitle    = {Abhandlungen aus dem mathematischen Seminar der Universit{\"a}t Hamburg},
  title        = {{\"U}ber ein Variationsprinzip zur L{\"o}sung von Dirichlet-Problemen bei Verwendung von Teilr{\"a}umen, die keinen Randbedingungen unterworfen sind},
  number       = {1},
  organization = {Springer},
  pages        = {9--15},
  volume       = {36},
  year         = {1971},
}

@Article{pavarino_2000,
  author    = {Pavarino, Luca F},
  title     = {Indefinite overlapping Schwarz methods for time-dependent Stokes problems},
  number    = {1-2},
  pages     = {35--51},
  volume    = {187},
  journal   = {Computer methods in applied mechanics and engineering},
  publisher = {Elsevier},
  year      = {2000},
}

@Article{embar_2010,
  author    = {Embar, Anand and Dolbow, John and Harari, Isaac},
  title     = {Imposing Dirichlet boundary conditions with Nitsche's method and spline-based finite elements},
  number    = {7},
  pages     = {877--898},
  volume    = {83},
  comment   = {They estimate the Nitsche stabilization parameter using local generalized eigenvalue problems. They address elliptic and fourth-order problems. This is directly related to Griebel and Schweitzer (2003), where a global generalized eigenvalue problem was used instead.},
  journal   = {International journal for numerical methods in engineering},
  publisher = {Wiley Online Library},
  year      = {2010},
}

@Article{wittum_1989,
  author    = {Wittum, Gabriel},
  title     = {Multi-grid methods for Stokes and Navier-Stokes equations},
  number    = {5},
  pages     = {543--563},
  volume    = {54},
  journal   = {Numerische Mathematik},
  publisher = {Springer},
  year      = {1989},
}

@Article{mittal_2005,
  author    = {Mittal, Rajat and Iaccarino, Gianluca},
  title     = {Immersed boundary methods},
  pages     = {239--261},
  volume    = {37},
  journal   = {Annu. Rev. Fluid Mech.},
  publisher = {Annual Reviews},
  year      = {2005},
}

@Article{glowinski_2007,
  author    = {Glowinski, R and Kuznetsov, Yu},
  title     = {Distributed Lagrange multipliers based on fictitious domain method for second order elliptic problems},
  number    = {8},
  pages     = {1498--1506},
  volume    = {196},
  journal   = {Computer Methods in Applied Mechanics and Engineering},
  publisher = {Elsevier},
  year      = {2007},
}

@Article{klawonn_1998_overlapping,
  author    = {Klawonn, Axel and Pavarino, Luca F},
  title     = {Overlapping Schwarz methods for mixed linear elasticity and Stokes problems},
  number    = {1-4},
  pages     = {233--245},
  volume    = {165},
  journal   = {Computer Methods in Applied Mechanics and Engineering},
  publisher = {Elsevier},
  year      = {1998},
}

@Article{li_2008,
  author    = {Li, Jian and He, Yinnian},
  title     = {A stabilized finite element method based on two local Gauss integrations for the Stokes equations},
  number    = {1},
  pages     = {58--65},
  volume    = {214},
  journal   = {Journal of Computational and Applied Mathematics},
  publisher = {Elsevier},
  year      = {2008},
}

@Article{burman_2010_ghost,
  author   = {Erik Burman},
  title    = {Ghost penalty},
  doi      = {https://doi.org/10.1016/j.crma.2010.10.006},
  issn     = {1631-073X},
  number   = {21},
  pages    = {1217 - 1220},
  url      = {http://www.sciencedirect.com/science/article/pii/S1631073X10002827},
  volume   = {348},
  journal  = {Comptes Rendus Mathematique},
  year     = {2010},
}

@Article{gander_2008,
  author  = {Gander, Martin Jakob},
  title   = {Schwarz methods over the course of time},
  pages   = {228--255},
  volume  = {31},
  journal = {Electronic transactions on numerical analysis},
  year    = {2008},
}

@Article{vanka_1986,
  author    = {Vanka, S Pratap},
  title     = {Block-implicit multigrid solution of Navier-Stokes equations in primitive variables},
  number    = {1},
  pages     = {138--158},
  volume    = {65},
  journal   = {Journal of Computational Physics},
  publisher = {Elsevier},
  year      = {1986},
}

@Article{deprenter_2019_multigrid,
  author    = {de Prenter, Frits and Verhoosel, Clemens V and van Brummelen, EH and Evans, JA and Messe, Christian and Benzaken, Joseph and Maute, Kurt},
  title     = {Multigrid solvers for immersed finite element methods and immersed isogeometric analysis},
  pages     = {1--32},
  volume    = {65},
  journal   = {Computational Mechanics},
  publisher = {Springer},
  year      = {2019},
}

@Article{flemisch_2007,
  author    = {Flemisch, Bernd and Wohlmuth, Barbara I},
  title     = {Stable Lagrange multipliers for quadrilateral meshes of curved interfaces in 3D},
  number    = {8},
  pages     = {1589--1602},
  volume    = {196},
  journal   = {Computer Methods in Applied Mechanics and Engineering},
  publisher = {Elsevier},
  year      = {2007},
}

@Article{duster_2008,
  author    = {D{\"u}ster, Alexander and Parvizian, Jamshid and Yang, Zhengxiong and Rank, Ernst},
  title     = {The finite cell method for three-dimensional problems of solid mechanics},
  number    = {45-48},
  pages     = {3768--3782},
  volume    = {197},
  journal   = {Computer methods in applied mechanics and engineering},
  publisher = {Elsevier},
  year      = {2008},
}

@Article{parvizian_2007,
  author    = {Parvizian, Jamshid and D{\"u}ster, Alexander and Rank, Ernst},
  title     = {Finite cell method},
  number    = {1},
  pages     = {121--133},
  volume    = {41},
  journal   = {Computational Mechanics},
  publisher = {Springer},
  year      = {2007},
}

\end{document}